\documentclass[conference]{IEEEtran}
\IEEEoverridecommandlockouts
\usepackage{cite}
\usepackage{amsmath,amssymb,amsfonts}
\usepackage{algorithmic}
\usepackage{graphicx}
\usepackage{textcomp}
\usepackage{xcolor}
\usepackage{listings}
\def\BibTeX{{\rm B\kern-.05em{\sc i\kern-.025em b}\kern-.08em
    T\kern-.1667em\lower.7ex\hbox{E}\kern-.125emX}}
\begin{document}

\title{Modelling Arbitrary Complex Dielectric Properties -- an automated implementation for gprMax}


\author{\IEEEauthorblockN{Sylwia Majchrowska\IEEEauthorrefmark{1},
Iraklis Giannakis\IEEEauthorrefmark{2},
Craig Warren\IEEEauthorrefmark{3}, and
Antonios Giannopoulos\IEEEauthorrefmark{4}}
\IEEEauthorblockA{\IEEEauthorrefmark{1}Department of Optics and Photonics,
Wroclaw University of Science and Technology, Wrocław, Poland}
\IEEEauthorblockA{\IEEEauthorrefmark{2}School of Geosciences, University of Aberdeen, Aberdeen, United Kingdom}
\IEEEauthorblockA{\IEEEauthorrefmark{3}Department of Mechanical and Construction Engineering,
Northumbria University,
Newcastle, United Kingdom}
\IEEEauthorblockA{\IEEEauthorrefmark{4}School of Engineering, The University of Edinburgh, Edinburgh, United Kingdom}
sylwia.majchrowska@pwr.edu.pl, iraklis.giannakis@abdn.ac.uk, craig.warren@northumbria.ac.uk, a.giannopoulos@ed.ac.uk
}

\maketitle

\begin{abstract}
There is a need to accurately simulate materials with complex electromagnetic properties when modelling Ground Penetrating Radar (GPR), as many objects encountered with GPR contain water, e.g. soils, curing concrete, and water-filled pipes. 
One of widely-used open-source software that simulates electromagnetic wave propagation is gprMax. It uses Yee's algorithm to solve Maxwell’s equations with the Finite-Difference Time-Domain (FDTD) method.
A significant drawback of the FDTD method is the limited ability to model materials with dispersive properties, currently narrowed to specific set of relaxation mechanisms, namely multi-Debye, Drude and Lorentz media. Consequently, modelling any arbitrary complex material should be done by approximating it as a combination of these functions.
This paper describes work carried out as part of the Google Summer of Code (GSoC) programme 2021 to develop a new module within gprMax that can be used to simulate complex dispersive materials using multi-Debye expansions in an automatic manner. The module is capable of modelling Havriliak-Negami, Cole-Cole, Cole-Davidson, Jonscher, Complex-Refractive Index Models, and indeed any arbitrary dispersive material with real and imaginary permittivity specified by the user.
\end{abstract}


\section{Introduction}
\label{sec:introduction}

gprMax (http://www.gprmax.com) was created in the mid-1990s as a tool that can be used to model Ground Penetrating Radar (GPR) responses from arbitrarily complex targets~\cite{bib:Giannopouls2005}.
gprMax was originally written in C programming language, and became one of the most widely used Finite-Difference Time-Domain (FDTD) based solvers in the GPR community.

In recent years gprMax has been re-developed and signiﬁcantly modernised to enable more complex models and advanced features~\cite{bib:Warren2016}.
At present, the open-source code is a combination of Python and Cython programming languages, which make it more flexible, user-friendly, and provide the opportunity to incorporate commonly used Python libraries such as Numpy and SciPy.
Moreover, the utilisation of Cython -- a superset of the Python language -- allows it to maintain a similar computation speed to the original C-based code.

The numerical modelling approach employed in gprMax is based on Yee’s~\cite{bib:Yee1966} algorithm (with second order accurate derivatives in space and time), which is used to solve Maxwell’s equations with the FDTD method.
However, the FDTD method is limited regarding modelling dispersive materials.
Due to its time-domain nature, arbitrary dispersive materials cannot be directly implemented in a straightforward manner~\cite{bib:Giannakis2014}.
To overcome this, the given dielectric spectrum is approximated by functions compatible with the time-domain numerical solver~\cite{bib:Kelley2007}.
The most common function employed for this is the multi-pole Debye expansion.
The package introduced in this work incorporates into gprMax three optimization approaches to fit a multi-pole Debye expansion to dielectric data.

The paper is organized as follows: Section~\ref{sec:materials} provides an overview of the properties of linear dispersive media, and the relaxation functions that were implemented; Section~\ref{sec:package} describes the design of the provided package along with the optimisation methods used; Section~\ref{sec:examples} presents the performance of the resulting Debye models in terms of computational time and dielectric model accuracy; and ﬁnally Section ~\ref{sec:summary} discusses the outcomes of the new implemented features.

\section{Dielectric Properties of Matter}
\label{sec:materials}

\subsection{Linear Dispersive Materials}

The dielectric properties of linear dispersive materials can be expressed as a frequency-dependent complex function 
\begin{equation}
\epsilon(\omega) = \epsilon^{'}(\omega) - j\epsilon^{''}(\omega),
\label{eq:permittivity}
\end{equation}
\begin{equation}
\mu(\omega) = \mu^{'}(\omega) - j\mu^{''}(\omega),
\label{eq:permeability}
\end{equation}
where $\omega$ is the angular frequency, $\epsilon^{'}$ and $\mu^{'}$ stand for the real parts of electric permittivity and magnetic permeability respectively; while $\epsilon^{''}$ and $\mu^{''}$ denote the imaginary parts.
In general, as a hard and fast rule, the real part dictates the velocity of the medium while the imaginary part is related to the electromagnetic losses~\cite{{bib:Balanis1989}}.
Values of complex relative permittivity and magnetic permeability can be obtained with a number of different techniques proposed and developed over the last decades~\cite{bib:afsar1986, bib:queffelec1994, bib:Baker-Jarvis1995}.

It appears that the most common group of materials, such as water, soil and biological tissues, exhibit dispersive behavior over a wide frequency range and linear behavior for commonly used field intensities~\cite{bib:taflove2005}.
The most important functions used to describe the behaviour of linear dispersive materials are the Debye, Lorentz, Drude, Havriliak-Negami~\cite{bib:Giannakis2014}, Jonsher equations, and the Complex Refractive Index Model (CRIM)~\cite{bib:Zadhoush2021}.

\subsection{Dielectric Relaxation Mechanisms}

\begin{figure}
\centerline{\includegraphics[width=\columnwidth]{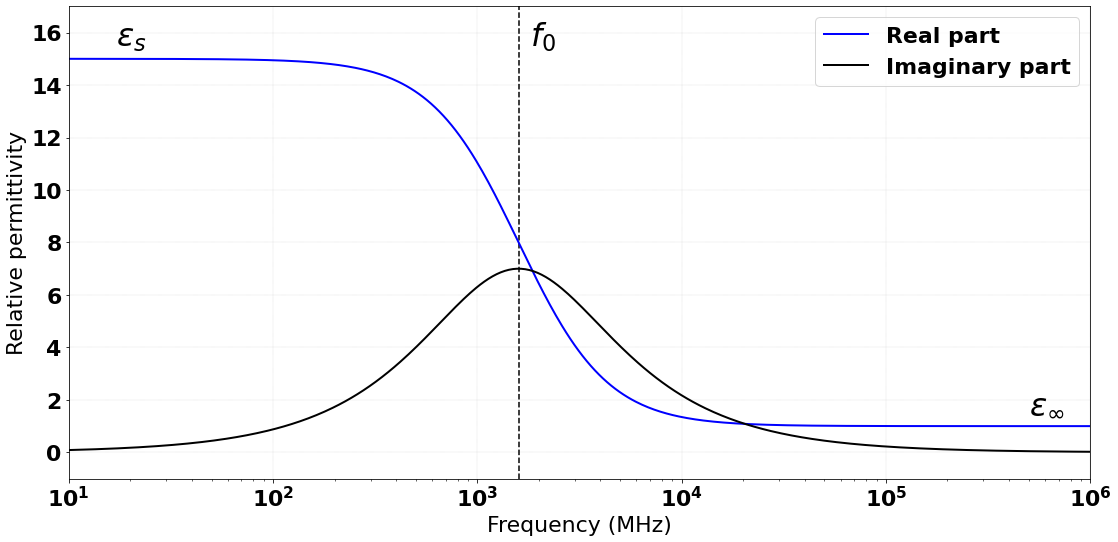}}
\caption{Both real and imaginary parts of a Debye function with \mbox{$\epsilon_{s} = 15\epsilon_{0}$} ($\epsilon_{0}$ is the electric permittivity of free space), \mbox{$\epsilon_{\infty} = \epsilon_{0}$} and \mbox{$\tau_{0} = 0.1\,\mathrm{ns}$}. The relaxation frequency equals \mbox{$f_{0} = \frac{1}{2\pi\tau_0}$}.}
\label{fig:single-Debye}
\end{figure}

Many materials of practical interest in GPR are linear, isotropic, and have insignificant magnetic response. In such cases they can be modeled using a scalar complex relative permittivity function.
Where frequency dependence of electric permittivity is required, typically the single-pole Debye equation is used \cite{bib:Balanis1989}
\begin{equation}
\epsilon(\omega) = \epsilon_{\infty} + \frac{\epsilon_{s} - \epsilon_{\infty}}{1+j\omega \tau_{0}},
\label{eq:Debye}
\end{equation}
where $\epsilon_{s}$ and $\epsilon_{\infty}$ are the values of the real part of the complex relative permittivity at static and infinite frequency, respectively, and $\tau_{0}$ is the relaxation time, which is directly related to relaxation frequency -- the frequency for which the imaginary part as well as the first derivative of the real part reach their maximum value (see Fig.~\ref{fig:single-Debye}).

However, the single-pole Debye equation is a poor model of dielectric behavior for most materials over wide frequency ranges.
To overcome this, the empirically derived Havriliak–Negami relaxation was proposed as a modification of the single-pole Debye relaxation model.
The equation additionally has two exponential parameters, and is expressed as
\begin{equation}
\epsilon(\omega) = \epsilon_{\infty} + \frac{\epsilon_{s} - \epsilon_{\infty}}{\left(1+\left(j\omega \tau_{0}\right)^{\alpha}\right)^{\beta}},
\label{eq:Havriliak–Negami}
\end{equation}
where $\alpha$ and $\beta$ are positive real constants ($0 \ge \alpha, \beta \le 1$). From this model, the Cole-Cole equation~\cite{bib:cole1941} setting $\beta = 1$, and Cole-Davidson equation~\cite{bib:manning1940} for $\alpha = 1$ can be derived. The Debye equation is obtained with $\alpha = 1$ and $\beta = 1$. The Havriliak-Negami function was first used to describe the dielectric properties of polymers~\cite{bib:havriliak_negami1967}, whereas the Cole-Cole equation is mainly used to model biological tissues~\cite{bib:ireland2013} and liquids~\cite{bib:barthel1995}.

On the other hand, the Jonscher function is mainly used to describe the dielectric properties of concrete~\cite{bib:bourdi2008} and soils\cite{bib:kjartansson1979}.
The frequency domain expression of the Jonscher function in the constant Q-factor approach (or quality Q-factor)~\cite{bib:bano2004} is given by
\begin{equation}
\epsilon(\omega) = \epsilon_{\infty} + A_{p}\left( -j\frac{\omega}{\omega_{p}} \right)^{n_p},
\label{eq:Jonsher}
\end{equation}
where $A_{p}$ is a Jonscher parameter (positive real constant), $n_{p}$ characterizes the change in amplitude as a function of frequency (it varies between 0 for materials with high dielectric loss and 1 for materials without dielectric losses), and $\omega_{p}$ is the reference frequency, arbitrarily chosen.

Other types of relaxation functions estimate the bulk permittivity of heterogeneous materials.
Here, CRIM~\cite{bib:birchak1974} was suggested as an improvement to the model initially proposed by Brown~\cite{bib:brown1956}.
The model illustrates the dielectric properties of the mixture with respect to the dielectric properties and the volumetric fractions of its components, and it is defined as
\begin{equation}
\epsilon(\omega)^{a} = \sum_{i=1}^{m}f_{i}\epsilon_{m,i}(\omega)^{a},
\label{eq:CRIM}
\end{equation}
where $a$ is a shape factor (commonly set to 0.5), $f_{i}$ and $\epsilon_{m,i}$ represents volumetric fraction and relative permittivity for \mbox{$i_\mathrm{th}$} material, respectively.
The CRIM has been established as the mainstream methodology in the GPR community mainly due to its simplicity and its straightforward implementation.

The Havriliak-Negami and Jonscher functions as well as steady Q-factor, CRIM, and experimentally obtained data cannot be directly implemented into the FDTD algorithm. An approach used by many researchers is the approximation of the complex permittivity using a multi-pole Debye expansion 
\begin{equation}
\epsilon(\omega) = \epsilon_{\infty} + \sum_{n=1}^N \frac{\Delta\epsilon_{n}}{1+j\omega \tau_{0, n}},
\label{eq:multi-Debye}
\end{equation}
with $\Delta\epsilon_{n}$ and $\tau_{0, n}$ being the change in permittivity and relaxation times of the Debye dispersion, respectively. Parameter $N$ indicates the number of Debye poles (the total number of Debye expansion components). 

\section{Package overview}
\label{sec:package}

\subsection{Code structure}

The current package features a stochastic global optimisation approach to fit a multi-pole Debye expansion to dielectric data.
The user can choose between Havriliak-Negami, Jonscher, and CRIM.
As mentioned previously, the Havriliak-Negami relaxation function is an inclusive function that holds as special cases the widely-used Cole-Cole and Cole-Davidson models.
The present package can be also used to fit arbitrary dielectric data derived experimentally or calculated using some formulated expression.

The package consists of two main modules:
\begin{itemize}
    \item \textit{Debye\_fit} -- being the main core of the package, containing the definition of the implemented relaxation functions and methods to run the optimization procedure,
    \item \textit{optimization} -- containing the implementation of selected global optimization approaches.
\end{itemize}

In the package we used three types of optimisation algorithms for determining the relaxation times, namely Dual Annealing (DA)~\cite{bib:tsallis1996}, Differential Evolution (DE)~\cite{bib:storn1997} and Particle Swarm Optimization (PSO)~\cite{bib:kennedy1995} techniques.
To calculate values of the weights we used the Damped Least-Squares (DLS) method. Our implementations are mainly based on modules provided in the pyswarm and SciPy optimize packages, which deliver various methods to minimize objective functions.

\subsection{Fitting procedure}

Accurately fitting Debye coefficients to a formulated relaxation function or a set of dielectric measurements across a very large frequency-range is problematic.
In literature various fitting methods have been investigated previously including genetic algorithms~\cite{bib:clegg2012, bib:krewer2013}, and hybrid particle-swarm least squares~\cite{bib:Kelley2007} optimization procedures. The whole idea is based on minimization of distance between the calculated values of the relaxation model and fitted multi-pole Debye expansion at individual frequency points in some defined range.

The optimisation method employed in the package is a hybrid linear-nonlinear optimisation approach. As a default we used a slightly adjusted optimization procedure to approximate the complex permittivity using the multi-pole Debye function expansion, as proposed by Kelley et. al.~\cite{bib:Kelley2007}. In their work, LS was used to determine the weights and PSO to determine the relaxation frequencies (\mbox{$\omega_{0,n} = 1/\tau_{0,n}$}).
Our modification is based on changing the sign to opposite in case of negative weights in order to introduce a large penalty in the optimisation process thus indirectly constraining the weights to always be positive.
Furthermore we added the real part to the cost action to avoid possible instabilities to arbitrary given functions that do not follow the Kramers–Kronig relationship.
These changes overcome some instability issues and thus make the process more robust and faster.

\begin{figure}
\centerline{\includegraphics[width=\columnwidth]{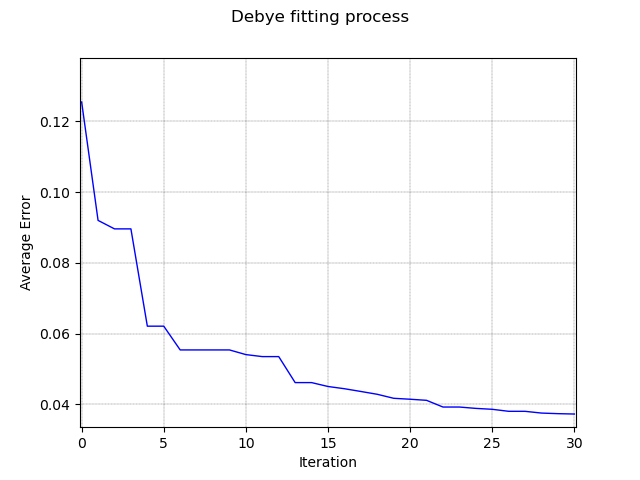}}
\caption{Example of fitting procedure for PSO algorithm.}
\label{fig:optimize}
\end{figure}

During creation of a relaxation object the user has a possibility to choose one from three available optimization techniques (and optionally some additional fixed parameters needed to completely specify the behavior of the optimizer) to calculate relaxation times. For the PSO procedure it is possible to plot the average error between actual and the approximated value of an objective function during the optimization process (see Fig.~\ref{fig:optimize}). Additionally we give a possibility to automatically set number of Debye poles. The estimation of the number of Debye poles is done iteratively starting from the one-pole Debye model.
The optimization stops when it reaches the desired accuracy (error below 5\%) or maximum number of iteration (20 Debye poles).
The procedure is launched when the parameter determining the number of Debye poles is set to -1.

\section{Multi-pole Debye models for FDTD simulations}
\label{sec:examples}

\begin{figure}
\centerline{\includegraphics[width=\columnwidth]{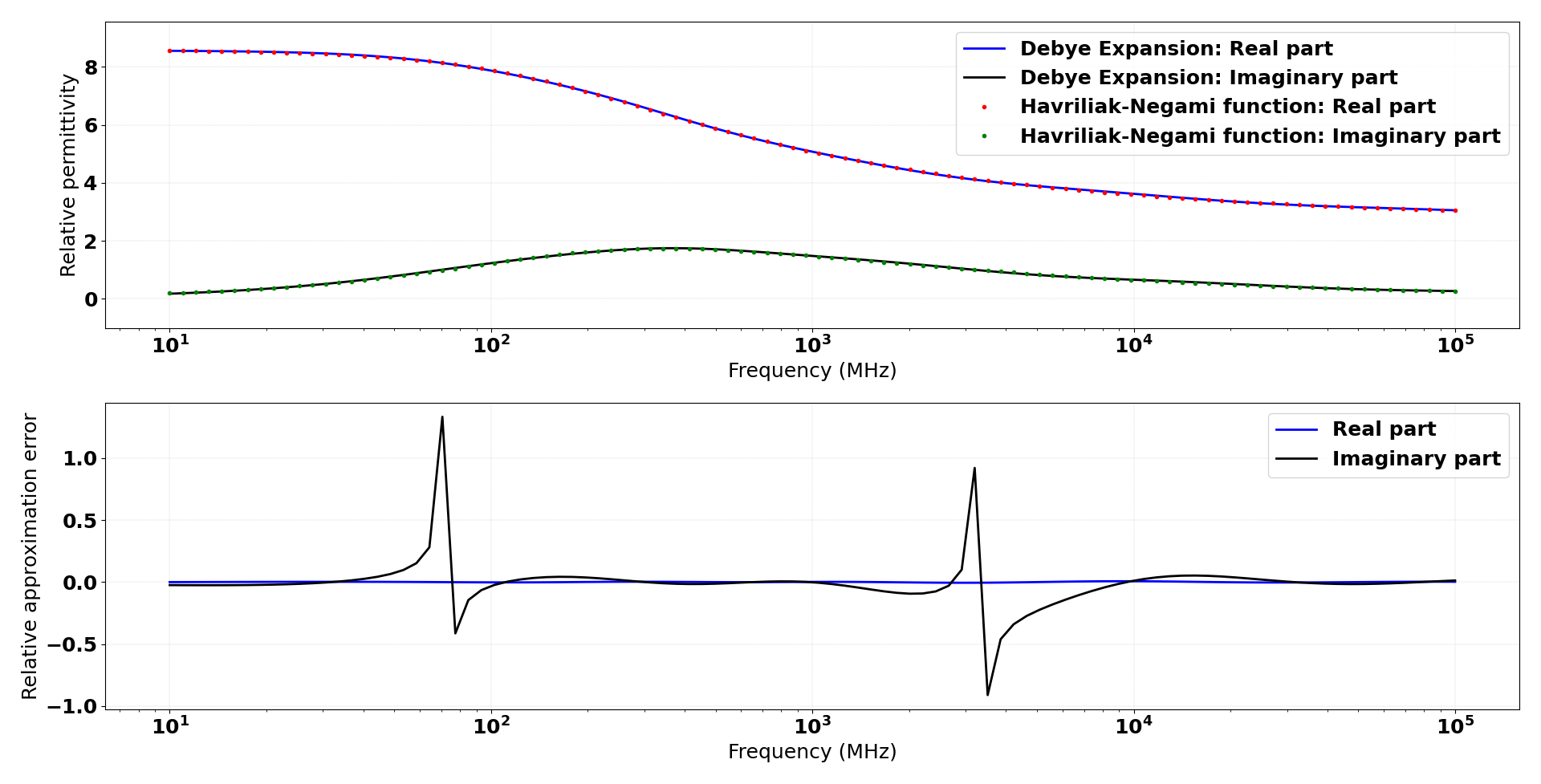}}
\caption{Real and imaginary parts of the five-pole Debye function expansions achieved with PSO algorithm compared to those of the Havriliak-Negami permittivity model~\cite{bib:Kelley2007}, and relative approximation error for both real and imaginary parts.}
\label{fig:HN}
\end{figure}

To illustrate the application of the package, a Havriliak-Negami permittivity model was approximated over the frequency range of 10 MHz to 100 GHz. The model parameters for the example are \mbox{$\epsilon_{s} = 8.6\epsilon_{0}$}, \mbox{$\epsilon_{\infty} = 2.7\epsilon_{0}$}, \mbox{$\tau_{0} = 0.94\,\mathrm{ns}$}, \mbox{$\alpha = 0.91$}, and \mbox{$\beta = 0.45$}. While the frequency range is wide, the number of Debye poles was set to 5 as in~\cite{bib:Kelley2007}. The results of the fitting procedure are presented in Fig.~\ref{fig:HN}.

We compared results from all algorithms in terms of computational speed and accuracy.
The average fractional error of each algorithm has been calculated separately for both real and imaginary part of the relative permittivity and then summed.
The obtained parameters are presented in Table~\ref{tab:fit}.
It appears that the PSO-DLS combination is the fastest optimization procedure, which achieves relatively good approximation results for the presented example (below 10\% of relative average error).

\begin{table}
\caption{Summary of results for fitting relative permittivity.}
\begin{center}
\begin{tabular}{|c|c|c|}
\hline
\textbf{Algorithm} & \textbf{Duration (s)} & \textbf{Average Error (\%)} \\ \hline
PSO-DLS & 0.97 &  8.48 \\
DA-DLS  & 2.56 &  4.37 \\
DE-DLS  & 3.61 &  4.48 \\
\hline
\end{tabular}
\label{tab:fit}
\end{center}
\end{table}

Multi-pole Debye, Drude, and Lorentz functions are already available in gprMax input files.
This made it easier to integrate this package into the software itself.
Listing~\ref{lst:HN} gives an example of the command to use a described above example using PSO method inside gprMax software.

\begin{lstlisting}[caption={A 5-pole Debye expansion to approximate Havriliak-Negami permittivity model.},label={lst:HN}]
#havriliak_negami: 1e7 1e11 0.91 0.45 2.7 5.9 9.4e-10 0.1 1 0 5 Kelley
\end{lstlisting}

The gprMax commands (like $\#havriliak\_negami$, $\#jonscher$, $\#crim$, $\#raw\_data$) to model dispersive materials deﬁne the basic material properties along with the relaxation function parameters. In Listing~\ref{lst:HN} the parameters for the $\#havriliak\_negami$ command are given going from left to right: the lower and upper frequency bounds (Hz), relaxation function parameters (here $\alpha$ and $\beta$ as for eq.~\ref{eq:Havriliak–Negami}), and basic material properties as real relative permittivity at infinite frequency, the difference between the static relative permittivity and the relative permittivity at infinite frequency, relaxation time (seconds), the conductivity (Siemens/metre), the relative permeability, the magnetic loss (Ohms/metre), the number of Debye poles, and an identifier for the material (text label). Additionally at the end of each command user can put optional integer, which controls the seeding of the random number generator used in stochastic global optimizer.
During fitting the hybrid PSO algorithm is used as a default.
Full description about other hashtag commands could be found in the documentation.

Listing~\ref{lst:output} gives an example of the output for the $\#havriliak\_negami$ command described above. Line 1 defines the basic material properties with the $\#material$ command, and in line 2 the $\#add\_dispersion\_debye$ command adds dispersive behaviour to the material based on the Debye formulation. In the definition of the $\#add\_dispersion\_debye$ command the number of Debye poles, and then the the difference between the static relative permittivity and the relative permittivity at infinite frequency for the $n_\mathrm{th}$ Debye pole, the relaxation time (seconds) for the $n_\mathrm{th}$ Debye pole are given. In case of running stand alone script these commands supposed to be added to the input file. On the other hand, the use of the $\#havriliak\_negami$ (Listing~\ref{lst:HN}) command will pass the approximated parameters inside gprMax automatically.

\begin{lstlisting}[caption={Example output for 5-pole Debye expansion to approximate Havriliak-Negami permittivity model.},label={lst:output}]
#material: 2.8345 0.1 1 0 Kelley
#add_dispersion_debye: 5 2.3563 4.3677e-10 0.6736 1.1623e-11 1.1009 1.5180e-9 1.2926 9.1048e-11 0.3091 1.1131e-12 Kelley
\end{lstlisting}

\section{Conclusions}
\label{sec:summary}

The current work describes the addition of a new advanced modelling feature in gprMax capable of simulating complex dielectric properties defined by some common functions such as Havriliak-Negami, Jonscher, CRIM, and any arbitrary user-specified dispersive dielectric properties. 
The hybrid linear-nonlinear optimisation approach employed here provides an effective, and accurate optimisation procedure to fit a multi-pole Debye expansion to the given dielectric data and subsequently implement it to gprMax in an automatic manner. 

\section*{Acknowledgment}

The project was funded via Google Summer of Code (GSoC) 2021 programme. GSoC initiative is a global program focused on bringing more student developers into open source software development.


\end{document}